\documentclass[11pt]{amsart}
\usepackage{ifthen,amssymb,graphicx}
\usepackage[all]{xy}
\usepackage{amsmath}
\usepackage{amsfonts}
\usepackage{mathrsfs}

  \newtheorem{thm}{Theorem}[section]
  \newtheorem{lem}[thm]{Lemma}
  
  \newtheorem{prop}[thm]{Proposition}
  
  \newtheorem{introthm}{Theorem}

  \theoremstyle{definition}
  \newtheorem{rem}[thm]{Remark}
  \newtheorem{defn}[thm]{Definition}
\newcommand{\dbc}[1]{{D}^b(#1)}
\newcommand{\what}[1]{{\widehat #1}}
\newcommand{\Hom}{{\operatorname{Hom}}}
\newcommand{\SHom}{{\mathcal{H}om}}
\DeclareMathOperator{\Pic}{{Pic}}
\newcommand{\cE}{{\mathcal E}}
\newcommand{\cF}{{\mathcal F}}
\newcommand{\cG}{{\mathcal G}}
\newcommand{\calL}{{\mathcal L}}
\newcommand{\cI}{{\mathcal I}}
\newcommand{\cO}{{\mathcal O}}
\newcommand{\cP}{{\mathcal P}}
\newcommand{\bR}{{\mathbf R}}
\newcommand{\cplx}[1]{{{\mathcal #1}^{\scriptscriptstyle\bullet}}}
\newcommand{\rk}{\operatorname{rk}}
\newcommand{\ch}{\operatorname{ch}}
\newcommand{\td}{\operatorname{td}}

\newcommand{\rest}[2]{{#1}_{\vert #2}}
\newcommand{\iso}{{\,\stackrel {\textstyle\sim}{\to}\,}}
\newcommand{\lra}{\longrightarrow}
\def\hra{\hookrightarrow}

\newcommand\op[1]{\operatorname{#1}}
\newcommand{\fm}{\mathscr S}
\title[CHARACTERIZATION OF JACOBIANS VIA PICARD BUNDLES]{A characterization of Jacobians by the existence of Picard Bundles}

\author[A. C. L\'opez Mart\'{\i}n]{Ana Cristina L\'opez Mart\'{\i}n}
\address{Departamento de Matem\'aticas {\rm and }Instituto Universitario de F\'{\i}sica Fundamental y Matem\'aticas
(IUFFYM), Universidad de Salamanca, Plaza
de la Merced 1-4, 37008 Salamanca, Spain}
\email{anacris@usal.es}
\email{dario@usal.es}

\author[E. C. Mistretta]{Ernesto Carlo Mistretta}
\address{Universit{\"a}t Bayreuth \\
Mathematisches Institut \\
Lehrstuhl Mathematik VIII \\
Universit{\"a}tstra{\ss}e 30 \\
D-95447 Bayreuth - Germany  }
\email{ernesto.mistretta@uni-bayreuth.de}

\author[D. S\'anchez G\'omez]{Dar\'{\i}o S\'anchez G\'omez}

\thanks{}
\date{\today}

\begin{document}
\begin{abstract}
Based on the Matsusaka-Ran criterion we give a criterion to characterize when a principal polarized abelian variety
is a Jacobian by the existence of Picard bundles.
\end{abstract}

\maketitle
%%%%%%%%%%

{\small \tableofcontents }

\section*{Introduction}
The problem of determining when an abelian variety is the Jacobian of a curve has been studied by many people
along the years. Generalizing the classical criterion of Matsusaka, Ran  gives in \cite{Ran81} a
characterization of Jacobians by the existence of curves with minimal cohomology class in the abelian variety.
This criterion is nowadays known as the Matsusaka-Ran criterion.

More recently, G. Pareschi and M. Popa use the theory of
Fourier-Mukai transforms as a useful tool in the study of the
existence of subvarieties of a principal polarized abelian variety
with minimal cohomology class. In this sense, they prove  in
\cite{PP08} a cohomological criterion which claims that if
$(A,\Theta)$ is an indecomposable principal polarized abelian
variety and $C$ is a geometrically non-degenerated reduced
equidimensional curve in $A$, such that the ideal sheaf
$\cI_C(\Theta)$ satisfies a Generic Vanishing property, then
$(A,\Theta)$ is the Jacobian of $C$ and $C$ has minimal cohomology
class. In the same paper they conjecture that if  the Index Theorem
with index 0 holds for $\cI_C(2\Theta)$,  with respect to the
Fourier-Mukai transform defined by the Poincar\'e bundle, then $C$
has minimal cohomology class. Consequently, using the Matsusaka-Ran
criterion, this would give a different  cohomogical criterion for
detecting Jacobians.

In this paper, we show the existence of such curves of minimal class using Picard bundles.

Picard bundles were introduced by Schwarzenberger in \cite{Schw63}
and have been used by many authors in the study of the geometry of
abelian varieties (c.f \cite{Kempf79a}, \cite{EL92}). Mukai studied
Picard bundles by means of  Fourier-Mukai transforms in (c.f.
\cite{Muk81}). We extend his definition of Picard bundles in order
to consider Fourier-Mukai transforms of any line bundle, and study
some properties of these sheaves in Proposition \ref{p:Picard}.

The fact that the locus of Jacobians is detected by the presence of Picard bundles appears in Kempf's work (cf.
\cite{Kempf79a} and \cite{Kempf82}). He shows that the projective bundle associated to a Picard bundle is the
symmetric product of the curve, and that deformations of Picard bundles are Picard bundles (this is proved
independently by Mukai). He shows furthermore that the deformations of Jacobians carrying a deformation of the
symmetric power of the curve (as a projective bundle on the Jacobian) are exactly those coming from deformations
of the curve. In this way, the presence of Picard bundles characterizes Jacobians among p.p.a.v.'s.

With the use of Fourier-Mukai transform and vanishing properties, we
carry out this description of the Jacobian locus explicitely,
determining those properties of Picard bundles that allow to detect
Jacobians. Moreover the Fourier-Mukai techniques give a quick and
clean way to recover known results about Picard bundles.

Our main result is the following:

% {\bf Theorem.} {\it Let $(A,\Theta)$ be an indecomposable p.p.a.v.
% of dimension $g$. If there exists a WIT$_g$ sheaf $\cF$ on $A$ with
% Chern classes $c_i(\cF)=(-1)^i\displaystyle\frac{\theta^i}{i!}$,
% then $(A,\Theta)$ is a Jacobian and a direct summand of $\cF$ is a
% Picard bundle.}

% \vspace{.5cm}

\begin{introthm}

Let $(A,\Theta)$ be an indecomposable p.p.a.v.
of dimension $g$. If there exists a WIT$_g$ sheaf $\cF$ on $A$ with
Chern classes $c_i(\cF)=(-1)^i\displaystyle\frac{\theta^i}{i!}$,
then $(A,\Theta)$ is a Jacobian and a direct summand of $\cF$ is a
Picard bundle.

\end{introthm}

\section*{Acknowledgements}  We are very grateful to G. Pareschi and M. Popa for proposing us this problem during
the summer school PRAGMATIC 2007. We would like to thank them for their lectures, useful conversations and the
constant encouragement. We would like to thank all the organizers of PRAGMATIC 2007 as well, for making possible
this very intense and interesting conference.

\section{Fourier-Mukai transforms for abelian varieties}
In this section, we recall some of the terminology of Fourier-Mukai
functors and the results that we will need in the rest of the paper.

Let $A$ be an abelian variety of dimension $g$ and $\what A=\Pic^0(A)$ the dual abelian variety.
$\what A$ represents the Picard functor, so there exists a universal
line bundle $\cP$ on $A\times \what A$, called the Poincar\'e bundle.
Thus, if $\alpha\in\what A$ corresponds to the  line bundle $\calL$ on $A$, one has
$${\cP}_\alpha:=\rest{{\cP}}{A\times\{\alpha\}}\simeq \calL\,.$$
Analogously, if $x\in A$, we denote
${\cP}_x:=\rest{{\cP}}{\{x\}\times\what A}\,$. The Poincar\'e bundle
can be normalised by the condition that
$\rest{{\cP}}{\{0\}\times\what A}$ is the trivial line bundle on
$\what A$. Denote $\pi_A\colon A \times \what A \to A$ and
$\pi_{\what A}\colon A\times \what A\to \what A$ the natural
projections.

The following result was proved by Mukai (cf. \cite{Muk81}).

\begin{thm}
The integral functor $\Phi\colon\dbc{A}\lra\dbc{\what A}$
defined by $\cP$
$$\Phi(\cplx{E}):=\bR{\pi_{\what A}}_\ast({\pi_{A}}^\ast(\cplx{E})\otimes\cP)$$
is a Fourier-Mukai transform, that is, an equivalence of categories.
Its quasi-inverse is the integral functor defined by ${\cP}^\ast[g]$
where ${\cP}^\ast$ denotes the dual of ${\cP}$.
\end{thm}

Let us denote by $\widehat{\Phi}\colon \dbc{\what A}\lra\dbc{A}$ the integral functor
defined by ${\cP}^\ast$. A straightforward consequence are the following isomorphisms:
$$\Phi\circ\what \Phi\simeq Id_{\dbc{\what{A}}}[-g]\quad\text{and}\quad \what \Phi\circ \Phi\simeq Id_{\dbc{A}}[-g]$$

\begin{rem}
In his original paper \cite{Muk81},
Mukai consider instead of $\widehat{\Phi}$ the integral functor $\fm\colon \dbc{\what A}\lra\dbc{A}$
defined by $\cP$.
The relation between the two  functors is given by $$\fm\simeq \widehat{\Phi}\circ (-1_{\what A})^\ast\, .$$
\end{rem}

For simplicity, we shall write $\Phi^j(\cplx{F})$ to denote the j-th cohomology sheaf of the complex $\Phi(\cplx{F})$,
and the same for the functor $\what \Phi$.
\begin{defn}
A coherent sheaf $\cF$ on $A$ is WIT$_i$ with respect to $\Phi$
(WIT$_i$-$\Phi$ for short) if $\Phi^j(\cF)=0$ for all $j\neq i$, or
equivalently if there exists a sheaf $\widehat{\cF}$ on $\what A$
such that $\Phi(\cF)\simeq\widehat{\cF}[-i]$. The sheaf
$\widehat{\cF}$ is called the Fourier-Mukai transform of $\cF$ with
respect to $\Phi$. When in addition $\widehat{\cF}$ is locally free,
we say that $\cF$ is IT$_i$ with respect to $\Phi$.
\end{defn}
We have analogous definitions of WIT and IT with respect the dual Fourier-Mukai functor $\what \Phi$.

The following proposition collects some properties of this special kind of sheaves.
\begin{prop}\label{prop:WIT&IT}
Let $\cF$ be a coherent sheaf on $A$. Then, the following holds:
\begin{itemize}
\item [1.] $\cF$ is IT$_i$-$\Phi$ if and only if $H^j(A,    \cF\otimes {\cP}_\alpha)=0$ for all $j\neq i$ and for all $\alpha\in\what A$.
\item [2.] $\cF$ is IT$_0$-$\Phi$ if and only if $\cF$ is WIT$_0$-$\Phi$.
\item [3.] If $\cF$ is WIT$_i$-$\Phi$, then $\widehat{{\cF}}$ is WIT$_{g-i}$-$\what \Phi$ and $\widehat{\widehat{{\cF}}\,}\simeq\cF$.
\item [4.] If $\cF$ is WIT$_g$-$\Phi$, then it is a locally free sheaf.
\item [5.] If $\cF$ is an ample line bundle, then it is IT$_0$-$\Phi$.
\end{itemize}
\end{prop}
\begin{proof}
Since $\cP$ is a locally free sheaf, $1)$ and $2)$ follow in a
straightforward way from Grauert's cohomology and base change theorem.
Part $3)$ follows from the isomorphism $\what \Phi\circ \Phi\simeq
[-g]$ and part $4)$  is a consequence of $ 3)$ and the definition of
IT. Part $5)$ is a consequence of the fact that ample line bundles
on abelian varieties have no higher cohomology (see for instance
\cite[III.16]{Mum74}) and 1).
\end{proof}
The relationship between the Chern characters of a WIT sheaf and those of its Fourier-Mukai transform is given by the following formula.

\noindent
\quad

\noindent
{\bf{Mukai's formula}}(\cite[Corollary 1.18]{Muk87b}):
If $\cE$ is a WIT$_j$-$\Phi$ sheaf, then
\begin{equation}\label{eq:Mukaiformula}
\ch_i(\widehat{\cE})=(-1)^{i+j}PD(\ch_{g-i}(\cE))
\end{equation}
where $PD$ denotes the Poincar\'e duality isomorphism.

\begin{defn}
A \textit{principally polarized abelian variety} (p.p.a.v. for short) is an abelian variety $A$ endowed with an
ample line bundle $\calL$ such that $\chi(\calL)=1$.
\end{defn}

\begin{rem}\label{r:ppav}
If $A$ is an abelian variety and we denote by $\tau_x$  the
translation morphism by a point $x\in X$, recall that $A$ is a
p.p.a.v. if and only if there exists an ample line bundle $\calL$ on
$A$ such that the morphism $\phi_{\calL}\colon A\to \what A$ defined
as $\phi_{\calL}(x)=\tau_{x}^\ast {\calL}\otimes {\calL}^{-1}$ is an
isomorphism. Moreover, by Proposition \ref{prop:WIT&IT}, the
polarization $\calL$ is IT$_0$, and it satisfies
\begin{equation}\label{eq:polarization}
\phi_{\calL}^\ast(\widehat{\calL})\simeq {\calL}^{-1}\, .
\end{equation}
\end{rem}

\section{Picard Bundles on Jacobians}
In this section we define Picard bundles via Fourier-Mukai transforms, and determine those properties of
theirs allowing us to characterize Jacobians among p.p.a.v.'s.

Let  $C$ be a smooth curve of genus $g\geq 2$ and consider $J_d(C)$ the Picard scheme parametrizing line bundles
of degree $d$ on $C$. This is a fine moduli space. Denote by $\mathcal{P}_d$ the universal Poincar\'{e} line
bundle on the direct product $C \times J_d(C)$ and $p\colon C\times J_d(C)\to C$ and $q\colon C\times J_d(C)\to
J_d(C)$ the projections. Fixing a point $x_0\in C$, it is normalized by imposing
$\rest{\mathcal{P}_d}{\{x_0\}\times J_d}\simeq \cO_{J_d}$. The higher direct images $\bR^i
q_{\ast}(\mathcal{P}_d)$ of $\mathcal{P}_d$ on $J_d(C)$ are known in the literature as degree $d$ Picard
sheaves.

Let us show how Picard sheaves can be seen in terms of the
Fourier-Mukai transform. Let $J_0(C)=J(C)$ be the Jacobian of $C$,
that is, the abelian variety that parametrizes the line bundles on
$C$ with degree zero.  The Riemann theta divisor $\Theta$ is a
natural polarization on $J(C)$ that defines a structure of
principally polarized abelian variety of dimension $g$ on $J(C)$. By
Remark \ref{r:ppav}, this gives a natural identification between
$J(C)$ and its dual abelian variety $\widehat{J(C)}$. With this
identification, if we denote by $$a\colon C\hookrightarrow J(C)$$
the Abel morphism, the normalized Poincar\'{e} bundle
$\mathcal{P}_0$ is precisely the restriction $(a\times 1)^\ast
\mathcal{P}$ of the universal line bundle $\mathcal{P}$ on
$J(C)\times J(C)$. On the other hand, the line bundle
$\mathcal{P}_d\otimes p^\ast\cO_C(-dx_0)$ defines an isomorphism
$\lambda_d\colon J_d(C)\iso J(C)$ and by normalization of the
Poincar\'{e} sheaves that we have considered, one has isomorphisms
$$\mathcal{P}_d \simeq (1\times \lambda_d)^\ast \mathcal{P}_0\otimes p^\ast \cO_C(dx_0)\, .$$

Using the base-change and the projection formulas, the Picard sheaf
$\bR^i q_{\ast}(\mathcal{P}_d)$ is $$\bR^i q_{
\ast}(\mathcal{P}_d)\simeq\lambda_d^\ast\Phi^i(a_\ast \cO_C(dx_0))\,
.$$

Considering that all Jacobians are already identified and although the last isomorphism is no longer true for an arbitrary line
bundle $L$ of degree $d$, the above discussion justifies the following
definition of Picard sheaves.
\begin{defn} Let $L$ be a line bundle on $C$ of degree $d$. The sheaves $\Phi^i(a_\ast L)$ are called the {\it degree $d$ Picard sheaves}.
\end{defn}

\begin{rem} The use of Fourier-Mukai transforms in the study of Picard bundles is originally due to Mukai \cite{Muk81}.
In this paper, he just considers the Picard sheaf $F_d=\Phi^1(a_\ast\cO_C(dx_0))$ corresponding to the line bundle $\cO_C(dx_0)$.
\end{rem}

\begin{rem}
Let $\Delta\colon \dbc{J(C)}\to \dbc{J(C)}$ be the dualizing functor
defined by  $\Delta(\cplx{F})=\bR \SHom(\cplx{F}, \cO_{J(C)})[g]$.
From Grothendieck duality, there is an isomorphism of functors
$$\Delta \circ\Phi\simeq ((-1)^\ast\circ \Phi\circ \Delta)[g]\, .$$
Taking into account that if $L$ is a line bundle on $C$, its derived
dual is $\Delta (a_\ast L)\simeq a_\ast(L^\ast \otimes \omega_C)[1]$
where $\omega_C$ is dualizing sheaf of $C$, one has an isomorphism
\begin{equation}\label{e:dual}\bR\SHom((\Phi(a_\ast L), \cO_{J(C))}))\simeq (-1)^\ast \Phi(a_\ast(L^\ast\otimes \omega_C))[1]\, .
\end{equation} which, in some cases, gives a duality relation between degree $d$ and degree $2g-2-d$  Picard bundles.
\end{rem}

Applying the theory of Fourier-Mukai transforms, we get some
properties of Picard sheaves that we summarize in the following
proposition
(cf. Theorem 4.2 and Proposition 4.3 in
\cite{Muk81}, properties of Picard sheaves as defined by Mukai).

\begin{prop}\label{p:Picard} The following holds:
\begin{enumerate}\item [1.] $\Phi^i(a_\ast L)$ are zero for $i\neq 0,1$.
\item [2.] For $d<0$, $\Phi^0(a_\ast L)=0$ and $\Phi^1(a_\ast L)$ is simple locally free of rank $g-d-1$. There is an isomorphism
$$\Phi^1(a_\ast L)\simeq (-1)^\ast \SHom((\Phi^0(a_\ast(L^\ast\otimes \omega_C)), \cO_{J(C)}))\, .$$
\item [3.] For $0\leq d<g-1$, $\Phi^1(a_\ast L)$ is supported on $J(C)$. %and its rank  at the generic point is $g-d-1$.
\item [4.] For $g-1\leq d<2g-1$, $\Phi^0(a_\ast L)$ and $\Phi^1(a_\ast L)$ are both non-zero.
\item [5.] For $d\geq 2g-1$,  $\Phi^0(a_\ast L)$ is a simple locally free sheaf of rank $d+1-g$ and $\Phi^1(a_\ast L)=0$.
There is an isomorphism
$$\Phi^0(a_\ast L)\simeq (-1)^\ast\SHom(\Phi^1(a_\ast(L^\ast\otimes \omega_C)), \cO_{J(C)})\, .$$
\end{enumerate}
\end{prop}

\begin{proof}
The first part is because the support of $L$ has dimension 1.
If $d<0$,
from Grauert's cohomology and base-change theorem  $\Phi^0(a_\ast L)=0$
and  $\Phi^1(a_\ast L)$ is a locally free sheaf of rank $g-d-1$.
 Since $\Phi$ is an equivalence of categories and $L$ is simple,
$\Phi(a_\ast L)[1]=\Phi^1(a_\ast L)$ is simple as well.
Analogously,
one gets the corresponding statements in 5).
In both cases, the duality relation between degree $d$ and degree $2g-2-d$ Picard bundles
is a consequence of the equation \eqref{e:dual}. Let us show 3). By cohomology and base-change,
one has that $\Phi^1(a_\ast L)_\alpha\simeq H^1(C, L\otimes \cP_\alpha)$ which,
 being $L\otimes \cP_\alpha$ of degree $d$,
is non-zero for every $\alpha \in J(C)$
because $\chi(L\otimes \cP_\alpha)<0$ by Riemman-Roch theorem. Now we prove 4).
By the equation \eqref{e:dual},
there is an isomorphism
$$\Phi^0(a_\ast L)\simeq (-1)^\ast \SHom^{-1}((\Phi(a_\ast(L^\ast\otimes \omega_C)), \cO_{J(C)}))\, .$$
From 3),
the sheaf $\Phi^1(a_\ast(L^\ast\otimes\omega_C))$ is supported on $J(C)$,
and  then the sheaf $\SHom^0(\Phi^{1}(a_\ast(L^\ast\otimes \omega_C)), \cO_{J(C)}))$ is non-zero.
 From the spectral sequence for local homomorphims
$$\begin{aligned}E_2^{p,q}=\SHom^p(\Phi^{-q}(a_\ast(L^\ast\otimes &\omega_C)),
\cO_{J(C)}))\Rightarrow \\ &E_{\infty}^{p+q}=\SHom^{p+q}(\Phi(a_\ast(L^\ast\otimes \omega_C)), \cO_{J(C)}))\, .
\end{aligned}$$
and the above isomorphism we easily conclude that $\Phi^0(a_\ast L)$ is non-zero.
Finally,
 to show that $\Phi^1(a_\ast L)$ is also non-zero,
consider the line bundle $L\otimes\cO_{C}(-dx_0)$,
where $x_0$ is the point of $C$ that we have fixed to normalize the Poincar\'{e} bundle.
This is a line bundle of degree zero and then $L\otimes\cO_{C}(-dx_0)\simeq \cP_\alpha$
 for some $\alpha\in J(C)$.
By Theorem 4.2 in \cite{Muk81},
there is a point $\kappa\in J(C)$ in the support of the sheaf $\Phi^1(\cO_C(dx_0))$.
Using again cohomology and base-change,
one obtains that $$H^1(C,\cO_C(dx_0)\otimes \cP_\kappa)\simeq H^1(C, L\otimes \cP_{\kappa-\alpha})$$
is non-zero. Hence, the point $\kappa-\alpha$ belongs to the support of $\Phi^1(a_\ast L)$
and we have the result.
\end{proof}

We will use the following
\begin{lem}\label{lem:c&ch}
Let $\cE$ be a vector bundle on a smooth variety $X$. The following are equivalent:
\begin{enumerate}
\item [a)] $\ch_j(\cE)=0\quad\text{ for all } j\geq 2$.
\item [b)] $\op{c}_i(\cE)=\displaystyle\frac{c_1(\cE)^i}{i!}\quad\text{ for all } i$.
\end{enumerate}
\end{lem}

\begin{proof}
By definition the total Chern class of $\cE$ is
$$\op{c}_t(\cE)=\op{c}_0(\cE)+\op{c}_1(\cE)t+\op{c}_2(\cE)t^2+\cdots+\op{c}_r(\cE)t^r=
\displaystyle\prod_{i=1}^r(1+a_it)$$
and the Chern character is
$$\ch(\cE)=\displaystyle\sum_{i=1}^r \displaystyle e^{a_it}=\sum\Big(1+a_it+\frac{(a_it)^2}{2!}+\cdots\Big)=r+\ch_1(\cE)t+\cdots$$
So, if $\op{c}_i(\cE)=\displaystyle\frac{c_1({\cE})^i}{i!}$ we get that $\op{c}_t(\cE)=\displaystyle e^{\op{c}_1(\cE)t}$. Then
\begin{align*}
\op{c}_1(\cE)t=& \log(\op{c}_t(\cE))=\sum \log(1+a_it)=\\
=&\sum (a_it-\frac{(a_it)^2}{2}+\frac{(a_it)^3}{3}+\cdots)=\\
=&\ch_1 (\cE)t-\ch_2(\cE)t^2+2\ch_3(\cE)t^3-3\ch_4(\cE)t^4+\cdots
\end{align*}
Hence $\ch_j(\cE)=0\quad\text{for any } j\geq 2$.

Conversely, if we assume that $\ch_j(\cE)=0\quad\text{for all } j\geq 2$, then one obtains
$$\log(\displaystyle\prod_{i=1}^r(1+a_it))=\op{c}_1(\cE)t$$ which implies
the condition b).
\end{proof}

Consider now the line bundle $L=\cO_C(2\Theta)\in J_{2g}(C)$.
By the proposition above,
the Picard sheaf $\Phi^0(a_\ast\cO_C(2\Theta))=\widehat {a_\ast\cO_C(2\Theta)}$
is a vector bundle on $J(C)$.

The aim of this section is to show some of the properties that this Picard bundle has:

\begin{prop}\label{thm:necessaryconditions}
If $\cF=\widehat {a_\ast\cO_C(2\Theta)}$, then
\begin{enumerate}
\item [1.] $\cF$ is a quotient of  $\widehat{\cO_{J(C)}(2\Theta)}$.
\item [2.] $\cF$ is WIT$_g$-$\widehat{\Phi}$.
\item [3.] $c_i(\cF)=(-1)^i\displaystyle\frac{\theta^i}{i!}$.
\item [4.] $\cF$ is simple.
\end{enumerate}
\end{prop}
\begin{proof}
Let us consider the exact sequence
\begin{equation}\label{eq:ideal}
0\lra\cI_C(2\Theta)\lra\cO_{J(C)}(2\Theta)\lra a_\ast\cO_C(2\Theta)\lra 0
\end{equation}
Since $\Theta$ is an ample divisor,  the line bundle
$\cO_{J(C)}(2\Theta)\otimes{\cP}_\alpha$ is also ample for any
$\alpha\in J(C)$. Thus, by applying the vanishing results for ample
line bundles on an abelian variety (see for instance
\cite[III.16]{Mum74}), we get that
$$H^i(J(C),\cO_{J(C)}(2\Theta)\otimes{\cP}_\alpha)=0\,\text{ for all  } i>0\text{ and all } \alpha\in J(C)$$
and, by Proposition \ref{prop:WIT&IT}, one concludes that $\cO_{J(C)}(2\Theta)$ is IT$_0$-$\Phi$.
On the other hand, Theorem 4.1 in \cite{PP03} proves that the sheaf $\cI_C(2\Theta)$ is IT$_0$-$\Phi$ as well.

By applying the Fourier-Mukai transform $\Phi$ to exact sequence (\ref{eq:ideal}) we then obtain
\begin{equation}\label{eq:idealFM}
0\lra\widehat{\cI_C(2\Theta)}\lra\widehat{\cO_{J(C)}(2\Theta)}\lra\widehat{a_\ast\cO_C(2\Theta)}\lra 0
\end{equation}
The next step is to compute the Chern classes of $\widehat{a_\ast\cO_C(2\Theta)}$.
In fact, this computation is originally due to Schwarzenberger \cite{Schw63}.
Here we deduce it using the Fourier-Mukai transform.

\noindent The Chern characters of $a_\ast\cO_C(2\Theta)$ can be
obtained using the Grothendieck Riemann Roch theorem for the Abel
morphism $C\stackrel{a}\hra J(C)$
$$\ch(a_\ast\cO_C(2\Theta))\cdot\td (J(C))=a_\ast(\ch(\cO_C(\rest{2\Theta}{C}))\cdot\td (C))$$
Remember that  $C$ has minimal cohomology class, that is,
$$[C]=\frac{\theta^{g-1}}{(g-1)!}$$ and the Todd class of $J(C)$ is trivial because it is an abelian variety. Then, one gets
\begin{equation}\label{eq:chernquotient}
\ch_j(a_\ast\cO_C(2\Theta))=
\begin{cases}
0 & j<g-1 \\
\displaystyle\frac{\theta^{g-1}}{(g-1)!} & j=g-1\\
g+1 & j=g
\end{cases}
\end{equation}
By applying Mukai's formula (\ref{eq:Mukaiformula}) we may compute the
Chern characters of Fourier-Mukai transform of
$a_\ast\cO_C(2\Theta)$. Thus we get
\begin{equation}\label{eq:chernquotientFM}
\ch_j(\widehat {a_\ast\cO_C(2\Theta)})=
\begin{cases}
g+1 & j=0\\
-\Theta & j=1\\
0 & j>1
\end{cases}
\end{equation}

Using Equation \eqref{eq:chernquotientFM} and Lemma \ref{lem:c&ch}, one obtains that
\begin{equation}
\op{c}_i(\widehat {a_\ast\cO_C(2\Theta)})=
\displaystyle(-1)^i\frac{\theta^i}{i!}
\end{equation}
Finally, notice that the sheaf $\cF$ is simple because
$$\Hom_{D(J(C))}(\widehat {a_\ast\cO_C(2\Theta)},\widehat
{a_\ast\cO_C(2\Theta)})\simeq \Hom_C(\cO_C(2\Theta),
\cO_C(2\Theta))\, ,$$ in particular, it is an indecomposable sheaf.
\end{proof}

\section{A Characterization of Jacobians via Picard Bundles}
In this section we shall use  the Matsusaka-Ran criterion to prove
that the existence of a sheaf satisfying properties 2) and 3) of the
Picard bundle in Proposition \ref{thm:necessaryconditions} is enough
to ensure that an indecomposable p.p.a.v. is the Jacobian of a
curve. Let us introduce some necessary notions and recall  the
Matsusaka-Ran criterion.

Let $C$ be an irreducible, or more generally connected, curve on an
abelian variety $A$. Some translate of $C$ generates
(group-theoretically) an abelian subvariety of $A$ which we denote
by $\langle C\rangle$ and we call the abelian subvariety generated
by $C$.
\begin{defn}\label{defn:generate}
A connected curve $C$ on an abelian variety $A$ is said to generate
$A$, if $\langle C\rangle=A$. More generally, an effective algebraic
1-cycle $\sum n_iC_i$ on $A$, with $n_i>0$ for all $n_i$, generates
$A$, if the union of the curves $C_i$ generates $A$.
\end{defn}
The following is the statement of the Matsusaka-Ran criterion (\cite[Theorem 11.8.1]{BirLa04}, see also
\cite[Theorem 3]{Ran81}).
\begin{thm}[Criterion of Matsusaka-Ran]\label{MatsusakaRan}
Suppose $(A,\Theta)$ is a polarized abelian variety of dimension $g$
and $C=\displaystyle\sum_{i=1}^r n_iC_i$ is an effective 1-cycle
generating $A$ with $[C]\cdot\Theta=g$. Then $n_i=1$ for all $1\leq
i\leq r$, the curves $C_i$ are smooth, and $(A,\Theta)$ is
isomorphic to the product of the canonically polarized Jacobians of
the $C_i's$:
$$(A,\Theta)\simeq (J(C_1),\Theta_1)\times\cdots\times (J(C_r),\Theta_r)$$
In particular, if $C$ is an irreducible curve generating $A$ with
$[C]\cdot \Theta=g$, then $C$ is smooth and $(A,\Theta)$ is the
Jacobian of $C$.
\end{thm}
Thus, the criterion that characterizes Jacobians by the existence of Picard bundles is the following
\begin{thm}\label{thm:sufficientconditions}
Let $(A,\Theta)$ be an indecomposable p.p.a.v. of dimension $g$. Suppose that there exists a sheaf $\cF$ on $A$
that satisfies the following conditions:
\begin{enumerate}
%\item [1.] $\sh{F}$ is a quotient of  $\cO_{\what A}(-2\what\Theta)$.
\item [1.] $\cF$ is WIT$_g$-$\widehat{\Phi}$.
\item [2.] $c_i(\cF)=(-1)^i\displaystyle\frac{\theta^i}{i!}$.
\end{enumerate}
Then there exists a smooth curve $C$ in $A$ such that
$(A,\Theta)\simeq (J(C),\Theta)$. Moreover,  if the sheaf $\cF$ is
indecomposable, then it is a simple degree $\rk(\cF)+g-1$ Picard
bundle with $\rk(\cF)\geq g$.
\end{thm}
\begin{proof} Consider $\widehat{\cF}$ the Fourier-Mukai transform of $\cF$.
Denote by $Z=\text{supp}(\widehat{\cF})$ the support of
$\widehat{\cF}$ and $i\colon Z\hookrightarrow A$ the natural
inclusion. Then, $\widehat{\cF}\simeq i_\ast \cG$ for some sheaf
$\cG$ on $Z$.

We recall that a WIT$_g$ sheaf is locally free (Proposition
\ref{prop:WIT&IT}). Using  Lemma \ref{lem:c&ch} and Mukai's formula
(\ref{eq:Mukaiformula}), we compute the components of the Chern
character of $\widehat{\cF}$ getting that
\begin{equation}
\ch_j(\widehat{\cF})=
\begin{cases}
0 & j<g-1 \\
\displaystyle\frac{\theta^{g-1}}{(g-1)!} & j=g-1\\
\rk(\cF) & j=g
\end{cases}
\end{equation}
This proves that $Z$ is a subscheme of codimension $g-1$. Define now
the 1-cycle $$Z_1(\widehat{\cF})=\sum_{\dim
V=1}l_V(\widehat{\cF})[V]$$ where the sum is over all 1-dimensional
subvarieties in $Z$ and $l_V(\widehat {\cF})$ is the length of the
stalk of $\widehat{\cF}$.

As a consequence of Grothendieck-Riemman-Roch theorem (c.f. \cite
[Theorem 18.3, Example 18.3.11]{Fu98}), it is  known that
$$\ch(\widehat{\cF})=Z_1(\widehat{\cF})+ \text{higher degree
terms}\, .$$ Hence the effective $1$-cycle $Z_1(\widehat{\cF})$ on
$A$ satisfies
$$Z_1(\widehat{\cF})=\displaystyle\frac{\theta^{g-1}}{(g-1)!}\,$$
This implies that this $1$-cycle generates $A$, by Corollary II.2,
Corollary II.3 and Lemma II.10 in \cite{Ran81}. Finally, since $(A,
\Theta)$ is indecomposable, and
$$[Z]=\displaystyle\frac{\theta^{g-1}}{(g-1)!}~,$$
then $Z$ is irreducible via the Poincar\'e duality. The Matsusaka-Ran
criterion allows us to conclude that the abelian variety
$(A,\Theta)$ is the Jacobian of a smooth curve, which proves the
first part of the theorem.

Assume now that the sheaf $\cF$ is indecomposable.
According to the previous discussion, the support $Z=C\sqcup W$,
where $C$ is the smooth curve and $W$ is a 0-dimensional closed subscheme.
Since $\widehat{\cF}\simeq i_\ast\cG$ is also indecomposable, then
$Z=C$, the inclusion
$i$ is $\pm a\colon C\hookrightarrow J(C)$, where $a$ is the Abel morphism,
and $\cG$ is a torsion free sheaf on $C$.
By applying Grothendieck-Riemman-Roch theorem to $\widehat{\cF}\simeq i_\ast\cG$,
we get
$$i_\ast (\rk(\cG))=\frac{\theta^{g-1}}{(g-1)!}\quad\text{ and } \quad i_\ast(c_1(\cG)-\frac{1}{2}\rk(\cG)K_C)=\rk(\cF)$$
$K_C$
being a canonical divisor of $C$.
Thus, $i=a$ and $\cG$ is a line bundle on $C$ of degree $\rk(\cF)+g-1$.
From  Proposition \ref{prop:WIT&IT} $a_\ast\cG$ is IT$_0$-$\Phi$ and then $\cF$ is simple and $\rk(\cF)\geq g$ by Proposition \ref{p:Picard}.
% This completes the proof.
\end{proof}

\begin{rem}
The same proof shows that when $(A,\Theta)$ is decomposable,
it is isomorphic to the direct product of the Jacobians
of the irreducible components of $Z_1(\widehat{\cF})$.
\end{rem}

\begin{rem}
When  $\cF$ is  not indecomposable,
we can show by the same argument that it is the direct sum of a Picard bundle
(the transform of  $\widehat{\cF}_{|C}$) and a vector bundle obtained as a chain of extensions of
degree $0$ line bundles (the transform of $\widehat{\cF}_{|W}$).
\end{rem}

\end{document}